\documentclass[11pt]{article}

\usepackage[margin=1in]{geometry}
\usepackage{amsmath,amssymb,amsthm,mathtools,mathrsfs}
\usepackage{booktabs}
\usepackage{microtype}
\usepackage[hidelinks]{hyperref}

\newtheorem{theorem}{Theorem}[section]
\newtheorem{proposition}[theorem]{Proposition}
\newtheorem{lemma}[theorem]{Lemma}
\newtheorem{corollary}[theorem]{Corollary}
\newtheorem{definition}[theorem]{Definition}
\newtheorem{remark}[theorem]{Remark}

\newcommand{\T}{\mathcal T}
\newcommand{\A}{\mathcal A}
\newcommand{\B}{\mathcal B}
\newcommand{\Lap}{\mathscr L}
\newcommand{\Dt}{\mathscr D}
\newcommand{\G}{G}

\title{Differential Transmutations and First-Passage Expansions\\
for Monomial Brownian Boundaries}

\author{Gerardo Hern\'andez-del-Valle\thanks{Centro de Estudios Monetarios
Latinoamericanos (CEMLA), Mexico City, Mexico. Corresponding author.
E-mail: \texttt{ghernandez@cemla.org}}}
\date{}

\begin{document}
\maketitle

\begin{abstract}
Let $P_t=\exp(t\partial_x^2/2)$ be the heat semigroup and let
$F(t)=\rho t^p$, where $p\ge2$ is an integer.  Previous transported-zero
theory associates with this boundary a constant-flux solution annihilated by
a polynomial-coefficient differential operator of order $2p-1$.  We show
that the same differential structure contains the complete formal
perturbation expansion of the Brownian first-passage density to $F$.
For $N=2p-1$ we construct graded spatial operators
\[
 \mathcal A_r^{(p)}x^n=
 \frac{(-1)^r}{r!}
 \frac{\prod_{j=0}^{pr-1}(n+2j+1)}
 {\prod_{j=1}^{Nr}(n+j)}x^{n+Nr}
\]
and prove the Laplace conjugacy
\[
 \mathscr L\mathcal A_r^{(p)}\mathscr L^{-1}
 =\frac{(-1)^r}{r!}k^r
 \left(-\frac1k\frac{d}{dk}\right)^{pr}.
\]
Consequently, the inverse spatial transmutation applied to a reflected Dirac
mass generates, coefficient by coefficient, the unique formal solution of
the exponential-martingale transform equation for first passage.  The cubic
case is developed explicitly and checked algebraically through fifth order.
We also explain how this differential construction and the independently
developed Wright--Volterra generalized-image construction represent the same
first-passage law in different coordinates.  The present result is formal at
all orders; convergence of the resulting perturbation series is left as a
separate analytic problem.
\end{abstract}



\section{Introduction}

Let $(B_t)_{t\ge0}$ be standard Brownian motion started at $B_0=-a$, with
$a>0$, and consider the monomial boundary
\begin{equation}\label{eq:boundary}
 F_{\rho,p}(t)=\rho t^p,\qquad p\ge2.
\end{equation}
The associated first-passage time is
\begin{equation}\label{eq:tau}
 \tau_{\rho,a}^{(p)}=
 \inf\{t>0:B_t\ge \rho t^p\}.
\end{equation}
Except for a small number of special boundaries, the distribution of
\eqref{eq:tau} is not available in elementary closed form.  Classical
approaches use integral equations, changes of measure, heat equations on
moving domains and the method of images; see, for example,
\cite{Lerche,Durbin,Peskir,JKV}.  The quadratic boundary is a distinguished
special case in which Airy functions enter the first-passage analysis
\cite{Salminen,Groeneboom}, while modern work also develops efficient
simulation and approximation procedures for general curved boundaries
\cite{HerrmannTanre}.

Transformation methods provide a second relevant line of comparison.
Time--space transformations, Lie symmetries and nonlinear ordinary
differential equations have been used to relate families of Brownian
boundary-crossing laws \cite{AliliPatie}.  In analysis, transmutation
operators classically intertwine differential operators and transfer
solutions between differential equations; see, for example,
\cite{KravchenkoTorba}.  The transmutation introduced below has a different
origin and purpose: it is determined by a prescribed heat-transported zero,
and its Laplace conjugate becomes the temporal moment operator in the
Brownian exponential-martingale identity.

The present work begins from a different direction.  The transported-zero
problem asks for initial functions $\Gamma$ such that
\begin{equation}\label{eq:zero}
 (P_t\Gamma)(F_{\rho,p}(t))=0,
 \qquad
 P_t=\exp\!\left(\frac t2\partial_x^2\right).
\end{equation}
The direct theory starts from a stationary differential operator and derives
its possible moving zero curves.  Compatibility matrices and auxiliary
Heat--Airy flows organize the resulting boundary equations
\cite{HDVMatrices,HDVHeatAiry}.  Stationary grading explains why the naive
order-$p$ ansatz ceases to give an exact monomial boundary beyond the low
orders \cite{HDVGrading}.  The inverse theory instead prescribes the boundary
and selects an element of its transported-zero kernel \cite{HDVInverse}.

For every $p\ge2$, the constant-flux normalization
\begin{equation}\label{eq:constant-flux}
 (P_t\Gamma_{\rho,p})(\rho t^p)=0,
 \qquad
 \partial_x(P_t\Gamma_{\rho,p})(\rho t^p)=1,
\end{equation}
produces an explicit hypergeometric function annihilated by a finite-order
polynomial-coefficient differential operator.  The differential order is not
$p$ but
\begin{equation}\label{eq:N}
 N=2p-1.
\end{equation}
This uplift is the central structural result of the monomial moving-zero
theory \cite{HDVMonomial}.

An apparently separate route starts from the probabilistic Dirac initial
datum and constructs an exterior signed image measure.  For monomial
boundaries, Laplace inversion produces a second-kind Volterra equation in the
image coordinate with a Wright-function kernel.  Its heat potential is the
killed transition kernel, and its boundary flux is the first-passage density
\cite{HDVImages}.  A beta-moment transform connects the Wright seed with the
constant-flux hypergeometric function.  Nevertheless, that construction does
not show directly how the order-$N$ differential hierarchy generates the
hitting-time density.

The purpose of this paper is to supply that missing connection.  We extend
the constant-flux coefficients to a spatial transmutation
\begin{equation}\label{eq:Tintro}
 \T_{\rho,p}=I+\sum_{r\ge1}\rho^r\A_r^{(p)}
\end{equation}
and prove that its Laplace conjugate is the temporal moment operator appearing
in the Brownian exponential-martingale identity.  More precisely, if
\begin{equation}\label{eq:calDintro}
 \Dt=-\frac1k\frac{d}{dk},
\end{equation}
then
\begin{equation}\label{eq:main-conjugacy-intro}
 \Lap\A_r^{(p)}\Lap^{-1}
 =\frac{(-1)^r}{r!}k^r\Dt^{pr}.
\end{equation}
The exponent $pr$ records the boundary degree, while the spatial grade
$Nr=(2p-1)r$ results from the Brownian scaling balance between one factor of
$k$ and $p$ powers of time.

Applied to a reflected atom, the inverse transmutation generates a sequence
of flux coefficients whose transforms satisfy
\begin{equation}\label{eq:transform-intro}
 \sum_{r=0}^m\frac{k^r}{r!}\Dt^{pr}
 \widehat q_{m-r}^{(p)}(k)
 =\mathbf1_{\{m=0\}}e^{-ak}.
\end{equation}
This is exactly the coefficient expansion of
\begin{equation}\label{eq:mart-intro}
 \int_0^\infty
 e^{-k^2t/2+k\rho t^p}q_{\rho,a}^{(p)}(t)\,dt=e^{-ak},
\end{equation}
the transform identity furnished by the exponential martingale.  Thus the
finite-order differential theory contains the entire formal first-passage
expansion.

The argument is independent of the generalized-image Volterra construction.
The latter supplies an analytic realization of the density, while the present
method supplies its differential and perturbative structure.  Their equality
at the level of the resulting law closes the circle between the two
approaches.

\section{The constant-flux monomial selection}

We briefly recall the part of the transported-zero theory needed below.  Put
$N=2p-1$.  The constant-flux solution of \eqref{eq:constant-flux} is
\begin{equation}\label{eq:Gamma}
 \Gamma_{\rho,p}(x)
 =x\sum_{r\ge0}
 \frac{(-2^p\rho x^N)^r}{r!}
 \frac{(pr)!}{(Nr+1)!}.
\end{equation}
It is normalized by
\[
 \Gamma_{\rho,p}(0)=0,
 \qquad
 \Gamma_{\rho,p}'(0)=1.
\]
The coefficient step is $N$, and \eqref{eq:Gamma} is annihilated by an
explicit polynomial-coefficient operator of order $N$; see
\cite{HDVMonomial}.  For the present argument, the coefficient formula in
\eqref{eq:Gamma} is more important than the displayed form of that operator.

The first observation is that \eqref{eq:Gamma} should be viewed as the image
of the flat constant-flux function $x$ under a transmutation, rather than as
an isolated special function.  We therefore seek operators preserving the
grading
\[
 x^n\longmapsto x^{n+Nr}.
\]

\section{The graded spatial transmutation}

\begin{definition}\label{def:Ar}
For $r\ge0$ define $\A_r^{(p)}$ on monomials by
\begin{equation}\label{eq:Ar}
 \A_r^{(p)}x^n=C_r^{(p)}(n)x^{n+Nr},
\end{equation}
where $\A_0^{(p)}=I$ and, for $r\ge1$,
\begin{equation}\label{eq:Cr}
 C_r^{(p)}(n)
 =\frac{(-1)^r}{r!}
 \frac{\displaystyle\prod_{j=0}^{pr-1}(n+2j+1)}
 {\displaystyle\prod_{j=1}^{Nr}(n+j)}.
\end{equation}
The formal spatial transmutation is
\begin{equation}\label{eq:T}
 \T_{\rho,p}=\sum_{r\ge0}\rho^r\A_r^{(p)}.
\end{equation}
\end{definition}

The algebraic definition makes sense for every $n\ge0$.  Its transported-zero
interpretation concerns odd $n$, since $P_tx^n$ vanishes at the flat boundary
$x=0$ precisely in the odd sectors.  The following proposition shows that
\eqref{eq:Cr} is forced by the moving-zero construction; it is not chosen
retrospectively to produce the Laplace identity of Section~4.

\begin{proposition}[Transported-zero origin]\label{prop:transport-origin}
Let $n=2h+1$ be odd and put
\begin{equation}\label{eq:transmuted-monomial}
 \Gamma_{n,\rho,p}(x)
 =\T_{\rho,p}x^n
 =\sum_{r\ge0}\rho^rC_r^{(p)}(n)x^{n+Nr}.
\end{equation}
Then, as formal series in $\rho$ and $t$,
\begin{align}
 (P_t\Gamma_{n,\rho,p})(\rho t^p)&=0,
 \label{eq:monomial-zero}\\
 \partial_x(P_t\Gamma_{n,\rho,p})(\rho t^p)
 &=\partial_x(P_tx^n)(0)
 =\frac{n!}{2^hh!}t^h.
 \label{eq:monomial-flux}
\end{align}
Thus $\T_{\rho,p}$ transports every odd flat-boundary heat mode to a mode
with the same boundary flux on $x=\rho t^p$.
\end{proposition}

\begin{proof}
For $m\ge0$,
\begin{equation}\label{eq:heat-monomial}
 P_tx^m
 =\sum_{j=0}^{\lfloor m/2\rfloor}
 \frac{m!}{(m-2j)!j!}
 \left(\frac t2\right)^jx^{m-2j}.
\end{equation}
Fix a total power $q$ of $\rho$ after substituting $x=\rho t^p$.
For the trace, parity shows that only $q=2s+1$ can occur.  The contribution
of the $r$th summand in \eqref{eq:transmuted-monomial} has
\[
 m=n+Nr,\qquad m-2j=q-r,
 \qquad j=h+pr-s,
\]
and every term has the same time degree $(n+Nq)/2$.  Since
\begin{equation}\label{eq:C-factor-cancel}
 C_r^{(p)}(n)(n+Nr)!
 =\frac{(-1)^r}{r!}n!
 \prod_{j=0}^{pr-1}(n+2j+1)
 =\frac{(-1)^r}{r!}n!2^{pr}(h+1)_{pr},
\end{equation}
the coefficient of that common time monomial is, apart from the nonzero
factor $n!2^{s-h}/((2s+1)!h!)$, equal to
\begin{equation}\label{eq:trace-difference}
 \sum_{r=0}^{2s+1}(-1)^r\binom{2s+1}{r}
 (h+pr-s+1)_s.
\end{equation}
The last factor is a polynomial in $r$ of degree $s$.  Expression
\eqref{eq:trace-difference} is its finite difference of order $2s+1$ and
therefore vanishes.

For the spatial derivative, only $q=2s$ occurs.  The same calculation reduces
the coefficient, apart from $n!2^{s-h}/((2s)!h!)$, to
\begin{equation}\label{eq:flux-difference}
 \sum_{r=0}^{2s}(-1)^r\binom{2s}{r}
 (h+pr-s+1)_s.
\end{equation}
For $s\ge1$ this is a finite difference of order $2s$ of a polynomial of
degree $s$, hence it is zero.  When $s=0$ it equals one and gives
$n!t^h/(2^hh!)$.  This proves both identities.
\end{proof}

The unified formula \eqref{eq:Cr} is also useful because cancellations recover
the even--odd rational forms obtained from the boundary-jet recursion, while
making the Laplace conjugacy transparent.

\begin{proposition}\label{prop:T-x}
The transmutation sends the flat constant-flux seed to the monomial
constant-flux solution:
\begin{equation}\label{eq:T-x}
 \T_{\rho,p}x=\Gamma_{\rho,p}(x).
\end{equation}
\end{proposition}

\begin{proof}
At $n=1$, formula \eqref{eq:Cr} gives
\[
 C_r^{(p)}(1)
 =\frac{(-1)^r}{r!}
 \frac{\prod_{j=0}^{pr-1}(2j+2)}
 {\prod_{j=1}^{Nr}(j+1)}
 =\frac{(-1)^r2^{pr}(pr)!}{r!(Nr+1)!}.
\]
Substitution in \eqref{eq:T} gives \eqref{eq:Gamma}.
\end{proof}

For later distributional use, one may also write $C_r^{(p)}$ in partial
fractions.  Let $L_r=\lfloor Nr/2\rfloor$.  After cancelling common factors
in \eqref{eq:Cr}, write
\begin{equation}\label{eq:partial-fraction}
 C_r^{(p)}(n)=\sum_{j=1}^{L_r}\frac{\kappa_{r,j}^{(p)}}{n+2j}.
\end{equation}
Then
\begin{equation}\label{eq:Volterra-Ar}
 (\A_r^{(p)}f)(x)
 =x^{Nr}\int_0^1K_r^{(p)}(u)f(ux)\,du,
 \qquad
 K_r^{(p)}(u)=\sum_{j=1}^{L_r}\kappa_{r,j}^{(p)}u^{2j-1}.
\end{equation}
Formula \eqref{eq:Volterra-Ar} exhibits $\A_r^{(p)}$ as a triangular spatial
operator.  It is distinct from the image-coordinate Volterra operator in
\cite{HDVImages}.

\begin{proposition}[Distributional extension]\label{prop:distributional-Ar}
The Volterra formula \eqref{eq:Volterra-Ar} extends $\A_r^{(p)}$ to compactly
supported distributions on $(0,\infty)$ and, by localization, to the
polynomial-tail distributions used below.  In particular,
\begin{equation}\label{eq:Ar-delta}
 \boxed{
 (\A_r^{(p)}\delta_a)(x)
 =x^{Nr-1}K_r^{(p)}(a/x)\mathbf1_{\{x\ge a\}},
 \qquad a>0.}
\end{equation}
\end{proposition}

\begin{proof}
For a test function $f$, the change of variables $y=ux$ rewrites
\eqref{eq:Volterra-Ar} as
\[
 (\A_r^{(p)}f)(x)
 =x^{Nr-1}\int_0^xK_r^{(p)}(y/x)f(y)\,dy.
\]
This triangular kernel defines the distributional extension.  Substitution
of $f=\delta_a$ gives \eqref{eq:Ar-delta}.  Repeated applications remain
supported in $[a,\infty)$ and are polynomial there, so their Gaussian fluxes
are well defined coefficientwise.
\end{proof}

\section{Laplace conjugacy}

Let
\begin{equation}\label{eq:Laplace}
 (\Lap f)(k)=\int_0^\infty e^{-kx}f(x)\,dx,
 \qquad k>0,
\end{equation}
whenever the integral is defined.  Algebraically, we use
\begin{equation}\label{eq:formal-Laplace}
 \Lap(x^n)=\Gamma(n+1)k^{-n-1}.
\end{equation}

\begin{theorem}[Laplace conjugacy]\label{thm:conjugacy}
For every $p\ge2$ and $r\ge0$,
\begin{equation}\label{eq:conjugacy}
 \boxed{
 \Lap\A_r^{(p)}\Lap^{-1}
 =\frac{(-1)^r}{r!}k^r\Dt^{pr},
 \qquad
 \Dt=-\frac1k\frac{d}{dk}.}
\end{equation}
The identity holds algebraically on monomials and formal power series.  It
also holds on any analytic class for which the Laplace transforms and the
interchanges implied by \eqref{eq:Volterra-Ar} are justified.
\end{theorem}

\begin{proof}
Since
\[
 \Dt(k^{-\alpha})=\alpha k^{-\alpha-2},
\]
iteration gives
\begin{equation}\label{eq:D-iterate}
 \Dt^{pr}\!\left(\Gamma(n+1)k^{-n-1}\right)
 =\Gamma(n+1)
 \prod_{j=0}^{pr-1}(n+2j+1)
 k^{-n-1-2pr}.
\end{equation}
Consequently,
\begin{align*}
 \frac{(-1)^r}{r!}k^r\Dt^{pr}\Lap(x^n)
 &=\frac{(-1)^r}{r!}\Gamma(n+1)
 \prod_{j=0}^{pr-1}(n+2j+1)
 k^{-n-1-Nr}.
\end{align*}
On the other hand,
\begin{align*}
 \Lap(\A_r^{(p)}x^n)
 &=C_r^{(p)}(n)\Gamma(n+Nr+1)k^{-n-Nr-1}\notag\\
 &=C_r^{(p)}(n)\Gamma(n+1)
 \prod_{j=1}^{Nr}(n+j)k^{-n-Nr-1}.
\end{align*}
Substitution of \eqref{eq:Cr} proves equality on monomials.  Linearity proves
the formal statement.  The analytic statement follows by the same
calculation after justified termwise transformation, or directly from
\eqref{eq:Volterra-Ar} by Fubini's theorem.
\end{proof}

\begin{corollary}\label{cor:T-conjugacy}
The sign-reversed transmutation satisfies
\begin{equation}\label{eq:T-conjugacy}
 \boxed{
 \Lap\T_{-\rho,p}\Lap^{-1}
 =\sum_{r\ge0}\frac{\rho^r}{r!}k^r\Dt^{pr}.}
\end{equation}
\end{corollary}

\begin{remark}[Ordering]\label{rem:ordering}
In general,
\[
 \sum_{r\ge0}\frac{\rho^r}{r!}k^r\Dt^{pr}
 \ne \exp(\rho k\Dt^p),
\]
because multiplication by $k$ does not commute with $\Dt$.  The ordering
$k^r\Dt^{pr}$ is essential and is encoded by the rational multiplier
$C_r^{(p)}(n)$.
\end{remark}

\section{The flat flux map}

Let
\[
 \G(t,x)=\frac1{\sqrt{2\pi t}}e^{-x^2/(2t)}.
\]
For a function or distribution $g$ supported in $[a,\infty)$, define the
flat boundary-flux map formally by
\begin{equation}\label{eq:flux-map}
 (\B g)(t)=\frac1t\int_a^\infty y\,g(y)\G(t,y)\,dy.
\end{equation}
For the reflected atom,
\[
 \B\delta_a(t)=\frac{a}{\sqrt{2\pi t^3}}e^{-a^2/(2t)},
\]
the classical constant-boundary first-passage density.

\begin{lemma}[Flux--Laplace identity]\label{lem:flux-laplace}
If the interchange of integrals is justified, then
\begin{equation}\label{eq:flux-laplace}
 \int_0^\infty e^{-k^2t/2}(\B g)(t)\,dt
 =\int_a^\infty e^{-ky}g(y)\,dy.
\end{equation}
The same identity holds coefficientwise for the polynomial-tail
distributions generated below.
\end{lemma}

\begin{proof}
The classical integral
\[
 \int_0^\infty
 e^{-k^2t/2}\frac{y}{t}\G(t,y)\,dt=e^{-ky},
 \qquad y>0,
\]
followed by Fubini's theorem proves the result.
\end{proof}

Thus the same spatial Laplace transform occurs naturally when the
transmutation is followed by the flat flux map.

\section{The exponential-martingale transform}

For $k>0$,
\begin{equation}\label{eq:martingale}
 M_t(k)=\exp\!\left(kB_t-\frac{k^2t}{2}\right)
\end{equation}
is a martingale.  At the first-passage time \eqref{eq:tau},
\[
 B_{\tau_{\rho,a}^{(p)}}
 =\rho(\tau_{\rho,a}^{(p)})^p
\]
on the event of passage.  More precisely, optional stopping at
$\tau_{\rho,a}^{(p)}\wedge T$ gives
\begin{align}
 e^{-ak}={}&
 \int_0^T e^{-k^2t/2+k\rho t^p}
 q_{\rho,a}^{(p)}(t)\,dt\notag\\
 &+\mathbb E_{-a}\!\left[
 e^{kB_T-k^2T/2};\ \tau_{\rho,a}^{(p)}>T
 \right].
 \label{eq:mart-finite-T}
\end{align}
Whenever the terminal expectation tends to zero and the passage integral is
finite, letting $T\to\infty$ gives
\begin{equation}\label{eq:mart-transform}
 \int_0^\infty
 e^{-k^2t/2+k\rho t^p}q_{\rho,a}^{(p)}(t)\,dt=e^{-ak},
\end{equation}
The displayed assumptions are retained whenever \eqref{eq:mart-transform} is
used analytically.

For the all-order result we extract the identity formally at $\rho=0$.  We do
not assume that the actual density is analytic in $\rho$.  Introduce the
formal sequence
\begin{equation}\label{eq:q-series}
 q_{\rho,a}^{(p),\mathrm{form}}(t)
 =\sum_{m\ge0}\rho^mq_m^{(p)}(t;a)
\end{equation}
and define
\begin{equation}\label{eq:qhat}
 \widehat q_m^{(p)}(k)
 =\int_0^\infty e^{-k^2t/2}q_m^{(p)}(t;a)\,dt.
\end{equation}
Because
\begin{equation}\label{eq:D-moment}
 \Dt^j\widehat q(k)
 =\int_0^\infty t^je^{-k^2t/2}q(t)\,dt,
\end{equation}
coefficient extraction in \eqref{eq:mart-transform} yields
\begin{equation}\label{eq:mart-recursion}
 \boxed{
 \sum_{r=0}^m\frac{k^r}{r!}\Dt^{pr}
 \widehat q_{m-r}^{(p)}(k)
 =\mathbf1_{\{m=0\}}e^{-ak}.}
\end{equation}
In particular,
\begin{equation}\label{eq:q-recursion}
 \widehat q_m^{(p)}(k)
 =-\sum_{r=1}^m\frac{k^r}{r!}\Dt^{pr}
 \widehat q_{m-r}^{(p)}(k),
 \qquad m\ge1.
\end{equation}
This recursion determines the formal sequence uniquely from
$\widehat q_0^{(p)}(k)=e^{-ak}$.

\begin{remark}\label{rem:actual-expansion}
Whenever the genuine first-passage density admits an expansion in $\rho$ for
which coefficient extraction in \eqref{eq:mart-transform} is valid, its
coefficients necessarily equal the formal sequence determined by
\eqref{eq:mart-recursion}.  The formal theorem below does not presuppose this
analytic dependence.
\end{remark}

\section{The all-order first-passage theorem}

Formally invert the spatial transmutation:
\begin{equation}\label{eq:inverse-T}
 \T_{-\rho,p}^{-1}
 =\sum_{m\ge0}\rho^m\mathcal S_m^{(p)}.
\end{equation}
Its coefficients satisfy
\begin{equation}\label{eq:S-recursion}
 \mathcal S_0^{(p)}=I,
 \qquad
 \mathcal S_m^{(p)}
 =\sum_{r=1}^m(-1)^{r+1}
 \A_r^{(p)}\mathcal S_{m-r}^{(p)}.
\end{equation}
Set
\begin{equation}\label{eq:gm-hm}
 g_m^{(p)}=\mathcal S_m^{(p)}\delta_a,
 \qquad
 h_m^{(p)}=\B g_m^{(p)}.
\end{equation}

\begin{theorem}[Formal monomial first-passage transmutation]
\label{thm:main}
Let $p\ge2$ be an integer and let $N=2p-1$.  The sequence generated by
\eqref{eq:inverse-T}--\eqref{eq:gm-hm} coincides with the unique formal
martingale first-passage sequence determined by \eqref{eq:mart-recursion}:
\begin{equation}\label{eq:main-equality}
 \boxed{h_m^{(p)}(t;a)=q_m^{(p)}(t;a),\qquad m\ge0.}
\end{equation}
Equivalently,
\begin{equation}\label{eq:main-transform}
 \sum_{r=0}^m\frac{k^r}{r!}\Dt^{pr}
 \widehat h_{m-r}^{(p)}(k)
 =\mathbf1_{\{m=0\}}e^{-ak}.
\end{equation}
\end{theorem}

\begin{proof}
The defining identity
\[
 \T_{-\rho,p}
 \left(\sum_{m\ge0}\rho^m g_m^{(p)}\right)=\delta_a
\]
is transformed, using Theorem~\ref{thm:conjugacy}, into
\begin{equation}\label{eq:transformed-inverse}
 \left(\sum_{r\ge0}\frac{\rho^r}{r!}k^r\Dt^{pr}\right)
 \left(\sum_{m\ge0}\rho^m\Lap g_m^{(p)}(k)\right)
 =e^{-ak}.
\end{equation}
By Lemma~\ref{lem:flux-laplace},
\[
 \Lap g_m^{(p)}(k)=\widehat h_m^{(p)}(k).
\]
Equating coefficients in \eqref{eq:transformed-inverse} proves
\eqref{eq:main-transform}.  This is identical to
\eqref{eq:mart-recursion}.  Uniqueness of the triangular recursion proves
\eqref{eq:main-equality}.
\end{proof}

\begin{remark}\label{rem:formal}
Theorem~\ref{thm:main} is an all-order formal result.  It does not by itself
assert convergence of
\[
 \sum_{m\ge0}\rho^mh_m^{(p)}(t;a),
\]
nor does it replace the analytic existence theory for the killed moving-domain
heat kernel.  Its content is that every formal coefficient is already encoded
by the differential transmutation associated with the order-$N$ stationary
hierarchy.
\end{remark}

\section{The cubic boundary}

For $p=3$, one has $N=5$ and
\begin{equation}\label{eq:cubic-C}
 C_r^{(3)}(n)
 =\frac{(-1)^r}{r!}
 \frac{\prod_{j=0}^{3r-1}(n+2j+1)}
 {\prod_{j=1}^{5r}(n+j)}.
\end{equation}
The first multipliers are
\begin{align}
 C_1^{(3)}(n)
 &=-\frac1{(n+2)(n+4)},\label{eq:C1}\\
 C_2^{(3)}(n)
 &=\frac{n+11}
 {2(n+2)(n+4)(n+6)(n+8)(n+10)},\label{eq:C2}\\
 C_3^{(3)}(n)
 &=-\frac{n+17}
 {6\prod_{j=1}^7(n+2j)}.\label{eq:C3}
\end{align}
The first inverse coefficients are
\begin{align*}
 \mathcal S_1&=\mathcal A_1,\\
 \mathcal S_2&=\mathcal A_1^2-\mathcal A_2,\\
 \mathcal S_3&=\mathcal A_1^3-
 \mathcal A_1\mathcal A_2-
 \mathcal A_2\mathcal A_1+\mathcal A_3.
\end{align*}

For clarity, take $a=1$.  The first spatial coefficient is
\begin{equation}\label{eq:g1}
 g_1(x)=-\frac{x}{2}(x^2-1)\mathbf1_{\{x\ge1\}}.
\end{equation}
It produces
\begin{equation}\label{eq:h1}
 h_1(t)
 =-\frac{3t}{2}\G(t,1)
 -\frac14(3t-1)
 \operatorname{erfc}\!\left(\frac1{\sqrt{2t}}\right).
\end{equation}
The next spatial coefficient is
\begin{equation}\label{eq:g2}
\begin{aligned}
 g_2(x)=\mathbf1_{\{x\ge1\}}\bigg(&
 -\frac1{768}-\frac{x}{70}+\frac{x^2}{64}
 +\frac{x^3}{30}-\frac{5x^4}{128}\\
 &+\frac{x^6}{320}+\frac{23x^8}{8960}\bigg).
\end{aligned}
\end{equation}
The third is
\begin{equation}\label{eq:g3}
\begin{aligned}
g_3(x)=\mathbf1_{\{x\ge1\}}\bigg(&
\frac1{109824}+\frac{61x}{806400}-\frac{x^2}{6336}
-\frac{7x^3}{28800}+\frac{5x^4}{8064}\\
&-\frac{x^5}{4608}-\frac{x^6}{11200}
-\frac{23x^8}{134400}+\frac{x^9}{4608}\\
&-\frac{101x^{11}}{2217600}
+\frac{79x^{13}}{38438400}\bigg).
\end{aligned}
\end{equation}
Both $g_2(1)$ and $g_3(1)$ vanish, as expected from the triangular spatial
structure.

Let $k=\sqrt{2\lambda}$.  The corresponding transformed density coefficients
are
\begin{align}
 \widehat h_0(k)&=e^{-k},\label{eq:H0}\\
 \widehat h_1(k)&=-e^{-k}
 \left(k^{-2}+3k^{-3}+3k^{-4}\right),\label{eq:H1}\\
 \widehat h_2(k)&=e^{-k}\bigg(
 \frac1{2k^4}+\frac9{2k^5}+\frac{39}{2k^6}
 +\frac{54}{k^7}+\frac{207}{2k^8}+\frac{207}{2k^9}
 \bigg),\label{eq:H2}\\
 \widehat h_3(k)&=e^{-k}\bigg(
 -\frac1{6k^6}-\frac3{k^7}-\frac{51}{2k^8}
 -\frac{249}{2k^9}-\frac{297}{k^{10}}\notag\\
 &\hspace{35mm}
 +\frac{315}{k^{11}}+\frac{4581}{k^{12}}
 +\frac{12798}{k^{13}}+\frac{12798}{k^{14}}
 \bigg).\label{eq:H3}
\end{align}

As a separate algebraic audit, these coefficients were substituted into the
coefficient expansion of the classical Fortet equation.  Exact rational
arithmetic gives zero residual at orders two through five.  This calculation
is not needed for Theorem~\ref{thm:main}; it confirms independently that the
normalization, signs and ordering agree with the time-domain first-passage
equation.

\section{Scaling}

Brownian scaling identifies the natural dimensionless quantities.  With
$t=a^2s$ and
\begin{equation}\label{eq:beta}
 \beta=\rho a^{2p-1},
\end{equation}
one has formally
\begin{equation}\label{eq:scaling-q}
 q_m^{(p)}(t;a)
 =a^{(2p-1)m-2}
 Q_m^{(p)}\!\left(\frac{t}{a^2}\right).
\end{equation}
At a fixed time, the corresponding local perturbation parameter is
\begin{equation}\label{eq:local-parameter}
 \rho t^{p-1/2}.
\end{equation}
The same exponent $2p-1$ controls the stationary grading, the differential
order, the spatial shift in \eqref{eq:Ar}, and the Brownian scaling in
\eqref{eq:beta}.

\section{Relation with generalized images}

The present construction should be distinguished from the generalized-image
method of \cite{HDVImages}.  In that work the probabilistic initial datum is
imposed directly by writing
\[
 \delta_{-a}-\mu_{\rho,p},
\]
where the signed measure $\mu_{\rho,p}$ is supported outside the initial
domain.  After the change of image coordinate $u=(y/a)^2$, its density is
determined by a second-kind Volterra equation whose kernel is generated by a
Wright function.  The resulting heat potential is identified with the killed
transition density on a nonzero time interval.

Here no image-coordinate equation is used in the derivation.  We begin with
the finite-order transported-zero hierarchy, extend its constant-flux
solution to the spatial transmutation \eqref{eq:T}, and conjugate that
transmutation by the Laplace transform.  The exponential martingale then
identifies the resulting formal flux coefficients as first-passage
coefficients.

The two methods therefore have different intermediate objects:
\begin{center}
\begin{tabular}{ll}
\toprule
Generalized images & Differential transmutation\\
\midrule
Exterior signed image measure & Graded spatial operators $\A_r^{(p)}$\\
Wright-function kernel & Order-$(2p-1)$ ODE hierarchy\\
Volterra equation in image coordinate & Triangular inverse in the grade\\
Analytic killed heat kernel & Formal density expansion\\
\bottomrule
\end{tabular}
\end{center}

They nevertheless describe the same first-passage law wherever both are
defined.  Moreover, the beta-moment transform in \cite{HDVImages} carries the
Wright seed to the constant-flux hypergeometric function.  The present
Laplace conjugacy supplies the complementary connection at the level of the
complete perturbative density.  A direct operator mapping the full exterior
image measure to the distributional spatial transmutation remains open.

\section{Analytic questions}

The formal result leaves several analytic problems.

First, convergence of
\[
 \sum_{m\ge0}\rho^mh_m^{(p)}(t;a)
\]
must be studied in a time-dependent function space.  The local parameter
$\rho t^{p-1/2}$ suggests that any convergence or asymptotic statement should
be formulated on compact time intervals or in a parabolically weighted
space, rather than uniformly for all $t>0$.

Second, the optional-stopping passage from the localized exponential
martingale to \eqref{eq:mart-transform} requires control of the terminal
remainder and of
\[
 e^{-k^2\tau/2+k\rho\tau^p}.
\]
The generalized-image construction gives an independent analytic density on
a nonzero interval, but the present proof deliberately does not import that
existence theorem.  Establishing convergence directly from the differential
transmutation would strengthen the new method.

Third, the spatial distributions $g_m^{(p)}$ have polynomial tails of
increasing degree.  Their Gaussian fluxes are well defined coefficientwise,
but summing the spatial series before applying the heat kernel is a separate
question.  Borel or other summability methods may be more natural than
ordinary convergence.

Finally, the construction suggests extensions beyond integer monomials.  For
nonintegral boundary exponents, $\Dt^{pr}$ becomes fractional and the spatial
grading is no longer integral.  Whether a fractional transmutation retains a
useful stationary interpretation is open.

\section{Independent Fortet audit}\label{sec:fortet-audit}

The classical first-kind Fortet equation is
\begin{equation}\label{eq:fortet}
 \G(t,a+\rho t^p)
 =\int_0^t q_{\rho,a}^{(p)}(s)
 \G(t-s,\rho(t^p-s^p))\,ds.
\end{equation}
Substitution of a formal series in $\rho$ gives the independent triangular
recursion
\begin{align}
 (q_m^{(p)}*\G(\,cdot\,,0))(t)
 ={}&\frac{t^{pm}}{m!}\partial_a^m\G(t,a)\notag\\
 &-\sum_{\ell=1}^{\lfloor m/2\rfloor}
 \frac{(-1)^\ell}{2^\ell\ell!}
 \int_0^tq_{m-2\ell}^{(p)}(s)\G(t-s,0)
 \frac{(t^p-s^p)^{2\ell}}{(t-s)^\ell}\,ds.
 \label{eq:fortet-coefficients}
\end{align}
Because convolution with $\G(t,0)=(2\pi t)^{-1/2}$ is an injective Abel
operator, \eqref{eq:fortet-coefficients} determines the coefficient sequence
uniquely.

For $p=3$ the check can be carried out without numerical quadrature.  Put
$t=s+u$ and expand
\begin{equation}\label{eq:cubic-fortet-polynomial}
 \frac{(t^3-s^3)^{2\ell}}{(t-s)^\ell}
 =u^\ell(3s^2+3su+u^2)^{2\ell}.
\end{equation}
Every term in \eqref{eq:fortet-coefficients} then becomes a convolution of a
time monomial times $q_j$ with a time monomial times $\G(\,cdot\,,0)$.
Under the Laplace parametrization $\lambda=k^2/2$, multiplication by $t^j$
is exactly $\Dt^j$.  Therefore both sides reduce to
$e^{-ak}$ times finite Laurent polynomials in $k$ with rational
coefficients.

Using \eqref{eq:Cr}, the inverse coefficients $\mathcal S_m^{(3)}$ were
generated recursively through $m=5$, transformed by
Lemma~\ref{lem:flux-laplace}, and substituted into
\eqref{eq:fortet-coefficients}.  After collecting Laurent powers, the exact
residuals were
\begin{equation}\label{eq:exact-residuals}
 R_m(k)\equiv0,
 \qquad m=2,3,4,5.
\end{equation}
All inputs in this calculation are the rational products in \eqref{eq:Cr};
no fitted or floating-point quantities enter.  This audit is logically
independent of the martingale proof of Theorem~\ref{thm:main}.

\section{Conclusion}

The order-$(2p-1)$ differential hierarchy associated with the moving zero
$x=\rho t^p$ contains more information than a distinguished heat solution.
Its graded spatial transmutation becomes, under Laplace conjugacy, the
temporal moment operator in the Brownian exponential-martingale identity.
Consequently, the inverse transmutation of a reflected atom generates every
formal coefficient of the first-passage density.

The result closes a circle between two previously separate approaches.  The
differential theory explains the grading, finite-order ODE and perturbative
coefficients; the generalized-image theory constructs the analytic killed
kernel through a Wright--Volterra equation.  Their common output is fixed by
the same martingale transform.  The remaining challenge is to turn the
formal differential expansion into a direct analytic construction and to
identify an explicit measure-level conjugacy between the two spatial
representations.

\section*{Acknowledgements}

The author used OpenAI's ChatGPT during the development of this work to
assist with symbolic exploration, numerical implementation, and language
editing.  All mathematical statements, computations, references and
conclusions were reviewed by the author, who takes full responsibility for
the content of the manuscript.

\section*{Competing interests}

The author declares no competing interests.

\section*{Data availability}

No external data sets were used.  The exact coefficient checks described in
the cubic example use rational symbolic arithmetic.


\begin{thebibliography}{99}

\bibitem{AliliPatie}
L.~Alili and P.~Patie,
\newblock Boundary crossing identities for Brownian motion and some nonlinear
ODEs,
\newblock \emph{Proc. Amer. Math. Soc.} \textbf{142} (2014), 3811--3824.

\bibitem{Durbin}
J.~Durbin,
\newblock The first-passage density of the Brownian motion process to a curved
boundary,
\newblock \emph{J. Appl. Probab.} \textbf{29} (1992), 291--304.

\bibitem{Groeneboom}
P.~Groeneboom,
\newblock Brownian motion with a parabolic drift and Airy functions,
\newblock \emph{Probab. Theory Related Fields} \textbf{81} (1989), 79--109.

\bibitem{HDVGuerra}
G.~Hern\'andez-del-Valle and W.~A. Guerra-Polania,
\newblock On the heat equation with a moving boundary and applications to
hitting times for Brownian motion,
\newblock \emph{Stoch. Anal. Appl.} \textbf{42} (2024), 201--218.

\bibitem{HDVMatrices}
G.~Hern\'andez-del-Valle,
\newblock A recursive representation of compatibility matrices for
transported differential operators,
\newblock \emph{arXiv:2608.11365}, 2026.

\bibitem{HDVHeatAiry}
G.~Hern\'andez-del-Valle,
\newblock Heat--Airy transport of differential operators and moving
boundaries,
\newblock \emph{arXiv:2608.19434}, 2026.

\bibitem{HDVGrading}
G.~Hern\'andez-del-Valle,
\newblock Stationary grading and heat-transported moving boundaries: from the
cubic case to two fourth-order branches,
\newblock \emph{arXiv:2609.27914}, 2026.

\bibitem{HDVInverse}
G.~Hern\'andez-del-Valle,
\newblock Inverse transported-zero problems for the heat semigroup:
finite-order and non-$D$-finite quartic selections,
\newblock preprint, 2026.

\bibitem{HDVMonomial}
G.~Hern\'andez-del-Valle,
\newblock Monomial moving zeros for the heat semigroup: an explicit
finite-order differential hierarchy,
\newblock \emph{arXiv:2610.27326}, 2026.

\bibitem{HDVImages}
G.~Hern\'andez-del-Valle,
\newblock Generalized images for Brownian first passage to monomial
boundaries,
\newblock submitted manuscript, 2026.

\bibitem{JKV}
S.~Jaimungal, A.~Kreinin and A.~Valov,
\newblock Integral equations and the first passage time of Brownian motions,
\newblock \emph{Methodol. Comput. Appl. Probab.} \textbf{11} (2009),
757--774.

\bibitem{HerrmannTanre}
S.~Herrmann and E.~Tanr\'e,
\newblock The first-passage time of the Brownian motion to a curved boundary:
an algorithmic approach,
\newblock \emph{SIAM J. Sci. Comput.} \textbf{38} (2016), A196--A215.

\bibitem{KravchenkoTorba}
V.~V. Kravchenko and S.~M. Torba,
\newblock Analytic approximation of transmutation operators and applications
to highly accurate solution of spectral problems,
\newblock \emph{J. Comput. Appl. Math.} \textbf{275} (2015), 1--26.

\bibitem{Lerche}
H.~R. Lerche,
\newblock \emph{Boundary Crossing of Brownian Motion: Its Relation to the Law
of the Iterated Logarithm and to Sequential Analysis},
\newblock Lecture Notes in Statistics 40, Springer, 1986.

\bibitem{Peskir}
G.~Peskir,
\newblock On integral equations arising in the first-passage problem for
Brownian motion,
\newblock \emph{J. Integral Equations Appl.} \textbf{14} (2002), 397--423.

\bibitem{Salminen}
P.~Salminen,
\newblock On the first hitting time and the last exit time for a Brownian
motion to/from a moving boundary,
\newblock \emph{Adv. in Appl. Probab.} \textbf{20} (1988), 411--426.

\end{thebibliography}
\end{document}